\documentclass[12pt,leqno]{amsart}
\usepackage{amssymb}
\usepackage{amsthm}
\usepackage{tikz}

\newtheorem{theorem}{Theorem}

\def\bZ{\mathbb Z}
\def\bF{\mathbb F}

\def\bG{\mathbb G}
\def\fp{\mathfrak p}
\def\fq{\mathfrak q}

\DeclareMathOperator{\Hom}{Hom}

\DeclareMathOperator{\Pic}{Pic}

\DeclareMathOperator{\Image}{Im}

\begin{document}
\title[Addendum] {The Spiegelungssatz for the Carlitz module ; an Addendum to : On a problem \`a la Kummer-Vandiver for function fields}
\author{Bruno Angl\`es \and Lenny Taelman}

\address{
Universit\'e de Caen, 
CNRS UMR 6139, 
Campus II, Boulevard Mar\'echal Juin, 
B.P. 5186, 
14032 Caen Cedex, France.
}
\email{bruno.angles@unicaen.fr}

\address{ 
Mathematisch Instituut,
Universiteit Leiden,
P.O. Box 9512,
2300~RA Leiden, The Netherlands.
}
\email{lenny@math.leidenuniv.nl}
\maketitle

\section{Introduction}
Let $p$ be an odd prime number. Let $\mu_p$ be the set of $p$th roots of unity in an algebraic closure of $\mathbb Q$ and ${\rm Cl}(\mathbb Q(\mu_p))$ be the ideal class group of $\mathbb Q(\mu_p)$. Then, using  class field theory and Kummer theory, one can show (see \cite[p. 188--191]{WAS}) that there exists a Galois-equivariant morphism:
$${\rm Hom}_{\mathbb Z}({\rm Cl}(\mathbb Q(\mu_p)), \mu_p) \rightarrow {\rm Cl}(\mathbb Q(\mu_p))[p],$$
such that its kernel is a  cyclic $\bF_p[{\rm Gal}(\mathbb Q(\mu_p)/\mathbb  Q)]$-module which is Galois-isomorphic to a subgroup of $\mathbb Z[\mu_p]^*/(\mathbb Z[\mu_p]^*)^p$. 
Leopoldt's Spiegelungssatz \cite[Theorem 10.9]{WAS} implies that also the cokernel of the above map is a cyclic
$\bF_p[{\rm Gal}(\mathbb Q(\mu_p)/\mathbb  Q)]$-module.
\par

In this  note we prove a kind of  analogue of Leopoldt's Spiegelungssatz for cyclotomic function fields.
Such a result is implicitly contained in \cite{TAE2} but not explicitly formulated.\par
The authors thank David Goss for helpful comments and suggestions.\par
\section{Notation}
We use the same notation as in \cite{TAE2}. So $k$ is a finite field of $q$ elements, $p$ is its characteristic, $A=k[T]$, $K=k(T)$ and $K_\infty = k((T^{-1}))$. 
Let 
\[
	\phi\colon A\rightarrow {\rm End}_{k} \bG_{a,A}, \,\, T \mapsto (x \mapsto Tx + x^q)
\]
be the Carlitz module.  If $R$ is an $A$-algebra we denote by $C(R)$ the $k$-vector space $R$ equipped by the $A$-module structure induced by $\phi$. 

We now fix a maximal ideal $\fp \subset A$ of degree $d$. Let $L/K$ be the splitting field of the $\fp$-torsion of $C$ and let $R \subset L$ be the integral closure of $A$ in $L$. We write $\Delta ={\rm Gal}(L/K)$. This is a cyclic group of order $q^d-1$. In fact, let $\Lambda$ be the module of $\fp$-torsion points in $C(L)$. We have $\Lambda \cong A/\fp$ non-canonically, and $\Delta$ acts on $\Lambda$ via the ``Teichm\"uller'' character $\omega\colon \Delta \to (A/\fp)^\times$, which is an isomorphism. We refer the reader to \cite[Ch. 7]{GOS}  or to \cite[Ch. 12]{ROS} for the basic properties of the abelian extension $L/K$.

Set $L_\infty := K_\infty \otimes_K L$. Then the Carliz exponential (see \cite{GOS} chapter 3), $\exp_C,$ induces an $A[\Delta]$-morphism $\exp_C\colon L_\infty \rightarrow C(L_\infty).$ The class module attached to $R$ and $\phi$ is defined by:
\[
	H(R) =\frac{C(L_\infty)}{\exp_C(L_\infty)+C(R)}.
\]
It is shown in \cite{TAE1} that the $A[\Delta]$-module $H(R)$ is finite.

\section{Spiegelungssatz}

As  a consequence  of the results obtained in \cite{TAE2} in order to prove an analogue of the Herbrand-Ribet theorem (using different methods,  this result has recently been refined in \cite{ANG&TAE2}),  we have the following analogue of Leopoldt's Spiegelungssatz.

\begin{theorem}
\label{theorem1} 
There is a natural $(A/\fp)[\Delta]$-morphism 
$${\rm Hom}_A( H(R), \Lambda )\rightarrow A/\fp \otimes_{\bF_p} (\Pic R)[p]$$
whose kernel and cokernel are cyclic $(A/\fp)[\Delta]$-modules.
\end{theorem}

A \emph{cyclic} module is a (possibly trivial) module generated by one element.

Note that the map in the theorem relates the dual of the $\fp$-part of the class \emph{module} to the $p$-part of the class \emph{group}!

\begin{proof}  We show how to obtain this from the results in \cite{TAE2}. \par

Let $\fq$ be the unique prime of $R$ above $\fp$. Combining Theorems 2 and 6 from \cite{TAE2} we find an exact sequence of $(A/\fp)[\Delta]$-modules
\begin{equation}\label{seq1}
	0 \to \Hom_A(H(R), \Lambda) \overset{f}{\to}
	A/\fp \otimes_k \Omega_{R}^{c=1} \to
	\Omega_R/\fq^{q^d} \Omega_R,
\end{equation}
where $\Omega_R$ is the module of K\"ahler differentials on $R$ over $k$ and $c$ is the ($k$-linear) $q$-Cartier operator. 
Kummer theory gives a short exact sequence of $(A/\fp)[\Delta]$-modules
\begin{equation}\label{seq2}
	0 \to A/\fp\otimes_{\bZ} R^\times \overset{g}{\to} 
	A/\fp \otimes_k \Omega_{R}^{c=1} \to
	A/\fp \otimes_{\bF_p} (\Pic R)[p] \to 0,
\end{equation}
see \cite {TAE2}, exact sequence (2) and  section 6.

Combing the exact sequences (\ref{seq1}) and (\ref{seq2}) we find an equivariant $A/\fp$-linear map
\[
	\alpha\colon \Hom_A(H(R), \Lambda) \to A/\fp \otimes_{\bF_p} (\Pic R)[p].
\]
Observe that  we have an isomorphism of $(A/\fp)[\Delta]$-modules :
$${\rm Ker}\,\alpha \cong \Image f\cap \Image g.$$
Thus the  kernel of $\alpha $  is a isomorphic to a submodule of the cyclic $(A/\fp)[\Delta]$-module $A/\fp\otimes_{\bZ} R^\times$, and therefore cyclic.

It remains to show that the cokernel is cyclic.  Observe that we have a $(A/\fp)[\Delta]$-isomorphism:
$${\rm Coker}\,  \alpha \cong \frac{A/\fp \otimes_k \Omega_{R}^{c=1}}{\Image f+\Image g}.$$
In particular, the cokernel of $\alpha$ is a quotient of the cokernel of $f$. Note that the cokernel of $f$ is a 
submodule of the kernel of the following map of $(A/\fp)[\Delta]$-modules
\begin{equation}\label{map1}
	1-c^d \colon \Omega_R/\fq^{q^d}\Omega_R \to \Omega_R/\fq\Omega_R.
\end{equation}
So it suffices to show that the kernel of the above map $1-c^d$ is a cyclic $(A/\fp)[\Delta]$-module.

The $(A/\fp)[\Delta]$-module $\Omega_R/\fq^{q^d}\Omega_R$ has a natural decreasing filtration
\[
	\Omega_R/\fq^{q^d}\Omega_R \supset \cdots
	\supset \fq^{i}\Omega_R/\fq^{q^d}\Omega_R \supset \cdots \supset 0.
\]
The intermediate quotient $\fq^i\Omega_R/\fq^{i+1}\Omega_R$ is generated by $\lambda^id\lambda$. It isomorphic with $A/\fp$ with $\Delta$ acting via $\omega^{i+1}$. In particular, we have
\[
	\dim_{A/\fp}  \frac{\Omega_R}{\fq^d\Omega_{R}} (\chi) = 
	\begin{cases}
		1 & \text{ if $\chi \neq \omega$ } \\
		2 & \text{ if $\chi = \omega$ }.
	\end{cases}
\]
Now consider the element $d\lambda$ of $(\Omega_R/\fq^d\Omega_R)(\omega)$. We have $c(d\lambda)=0$
and $d\lambda$ is nonzero in $\Omega_R/\fq\Omega_R$. So we find that $(c^d-1)(d\lambda)$ is nonzero in $\Omega_R/\fq\Omega_R$ and conclude that the $\omega$-part of the kernel of the map (\ref{map1}) is one-dimensional. In particular, we find that the kernel of (\ref{map1}) is a cyclic $(A/\fp)[\Delta]$-module, and hence that the cokernel of $\alpha$ is cyclic.
\end{proof}

Let $\chi$ be a character $\Delta \to (A/\fp)^\times$ with $\chi(k^\times) \neq \{1\}$. Assume that
$(A/\fp \otimes_A H(R_P)) (\omega \chi^{-1})=\{0\}.$ The above theorem implies that
$(A/\fp \otimes_\bZ Cl^0(L) )(\chi)$ is a cyclic $A/\fp$-module.  This is a key ingredient in the construction \cite{ANG&TAE}  of counterexamples to the analogue of the Kummer-Vandiver conjecture for the class module. In \cite{ANG&TAE2} it is shown that these also yield counterexamples to Anderson's conjecture \cite{AND}.


  \end{document}